\documentstyle[twoside]{article}
\title{\bf An expansion for the sum of a product of an exponential and a Bessel function. II}
\author{\sc R.B. Paris\footnote{E-mail address:\ \ {\tt r.paris@abertay.ac.uk}}\\
\\
{\em Division of Computing and Mathematics,}\\
{\em Abertay University, Dundee DD1 1HG, UK}\\
}

\begin{document}
\newcommand{\bee}{\begin{equation}}
\newcommand{\ee}{\end{equation}}
\def\f#1#2{\mbox{${\textstyle \frac{#1}{#2}}$}}
\def\dfrac#1#2{\displaystyle{\frac{#1}{#2}}}
\newcommand{\fr}{\frac{1}{2}}
\newcommand{\fs}{\f{1}{2}}
\newcommand{\g}{\Gamma}
\newcommand{\br}{\biggr}
\newcommand{\bl}{\biggl}
\newcommand{\ra}{\rightarrow}
\renewcommand{\topfraction}{0.9}
\renewcommand{\bottomfraction}{0.9}
\renewcommand{\textfraction}{0.05}
\newcommand{\mcol}{\multicolumn}
\date{}
\maketitle
\pagestyle{myheadings}
\markboth{\hfill {\it R.B. Paris} \hfill}
{\hfill {\it A Bessel function sum. II } \hfill}
\begin{abstract} 
We examine the sum of a decaying exponential (depending non-linearly on the summation index) and a Bessel function in the form
\[\sum_{n=1}^\infty e^{-an^p}\frac{J_\nu(an^px)}{(\fs an^px)^\nu}\qquad (x>0),\]
where $J_\nu(z)$ is the Bessel function of the first kind of real order $\nu$ and $a$ and $p$ are positive parameters. By means of a Mellin transform approach we obtain an asymptotic expansion that enables the evaluation of this sum in the limit $a\to 0$. A similar result is derived for the sum when the Bessel function is replaced by the modified Bessel function $I_\nu(z)$ when $x\in (0,1)$. The case of even $p$ is of interest since the expansion becomes exponentially small in character. We demonstrate that in the case $p=2$, a result analogous to the Poisson-Jacobi transformation exists for the above sum.
\vspace{0.4cm}

\vspace{0.3cm}

\end{abstract}

\vspace{0.5cm}

\begin{center}
{\bf 1.  Introduction}
\end{center}
\setcounter{section}{1}
\setcounter{equation}{0}
\renewcommand{\theequation}{\arabic{section}.\arabic{equation}}
We consider the sum
\begin{equation}\label{e11}
S_{\nu,p}(a,x)=\sum_{n=1}^\infty e^{-an^p}\,\frac{J_\nu(an^px)}{(\fs an^px)^\nu},\qquad p>0,\ x>0,\ a>0,
\end{equation}
where $J_\nu(x)$ is the Bessel function of the first kind of real order $\nu$ and
\[(\fs z)^{-\nu}\!J_\nu(z)=\sum_{k=0}^\infty \frac{(-z^2/4)^k}{\g(1+\nu+k) k!}.\]
Our interest herein is the asymptotic expansion of $S_{\nu,p}(a,x)$ in the limit $a\to0$ when convergence of the above sum becomes slow. We employ a Mellin transform approach to express the sum as a Mellin-Barnes integral involving the Riemann zeta function $\zeta(s)$ and the Gauss hypergeometric function ${}_2F_1$. It is found that the resulting (algebraic) asymptotic series converges when $p<1$ and diverges when $p>1$; the case $p=1$ requires the condition
$a<2\pi/\sqrt{1+x^2}$ for convergence.

In \cite{P17} the  asymptotic expansion of the sum
\[\sum_{n=1}^\infty e^{-an} \frac{J_\nu(bn)}{(\fs bn)^\nu}\]
was examined for $a\to 0$ with $0<b<2\pi$ and $\nu>-\fs$. The details in this case relied on the use of a double Mellin-Barnes integral involving the zeta function when $a<b$. 
When $a=0$, the sum has been considered by Tri\u ckovi\'c {\it et al.} in \cite{TVS}, where approaches using Poisson's summation formula and Bessel's integral were employed to derive convergent expansions. 

An interesting situation arises when $p$ is an even integer in the sum (\ref{e11}), where the algebraic asymptotic series vanishes to leave an exponentially small contribution. We pay particular attention to the case $p=2$, where it will be demonstrated that a transformation analogous to the well-known Poisson-Jacobi transformation \cite[p.~124]{WW}
\bee\label{e12}
\sum_{n=1}^\infty e^{-an^2}=\frac{1}{2}\sqrt{\frac{\pi}{a}}-\frac{1}{2}+\sqrt{\frac{\pi}{a}} \sum_{n=1}^\infty e^{-\pi^2 n^2/a} \qquad (\Re (a)>0)
\ee 
holds for the sum in (\ref{e11}).

\vspace{0.6cm}

\begin{center}
{\bf 2.  The asymptotic expansion of $S_{\nu,p}(a,x)$ for $a\to0$}
\end{center}
\setcounter{section}{2}
\setcounter{equation}{0}
\renewcommand{\theequation}{\arabic{section}.\arabic{equation}}
Let $a>0$, $x>0$, $p>0$ and $\nu$ be a real parameter. We consider the asymptotic expansion of the sum
\bee\label{e21}
S_{\nu,p}(a,x)=\sum_{n=1}^\infty e^{-an^p} \,\frac{J_\nu(an^px)}{(\fs an^px)^\nu}
\ee
for $a\to0$. We adopt a Mellin transform approach as described, for example, in \cite[p.~118]{PK} and write the above sum as
\[S_{\nu,p}(a,x)=(\fs x)^{-\nu} \sum_{n=1}^\infty h(an^p),\qquad h(\tau):=e^{-\tau} \tau^{-\nu} J_\nu(x\tau).\]
With the Mellin transform given by $H(s)=\int_0^\infty \tau^{s-1}h(\tau)\,d\tau$, we have \cite[p.~385(2)]{WBF}
\[H(s)=\int_0^\infty \tau^{s-\nu-1}e^{-\tau} J_\nu(x\tau)\,d\tau
= \frac{(\fs x)^\nu\g(s)}{\g(1+\nu)}\,{}_2F_1(\fs s,\fs s+\fs;1+\nu;-x^2),\]
where ${}_2F_1(\alpha,\beta;\gamma;z)$ is the Gauss hypergeometric function \cite[p.~384]{DLMF}. Then upon use of the Mellin inversion theorem \cite[p.~118]{PK} we obtain the integral representation in the form
\[S_{\nu,p}(a,x)=(\fs x)^{-\nu} \sum_{n=1}^\infty \frac{1}{2\pi i}\int_{c-\infty i}^{c+\infty i} H(s) (an^p)^{-s}ds=\frac{(\fs x)^{-\nu}}{2\pi i} \int_{c-\infty i}^{c+\infty i} H(s)\,\zeta(sp)\,a^{-s}ds\]
\bee\label{e22}
=\frac{1}{\g(1+\nu)}\,\frac{1}{2\pi i}\int_{c-\infty i}^{c+\infty i} \g(s)\,\zeta(sp) \,{}_2F_1(\fs s,\fs s+\fs;1+\nu;-x^2)\,a^{-s}ds,
\ee
where $\zeta(s)$ is the Riemann zeta function and the integration path is such that $c>1/p$.

The integrand in (\ref{e22}) has simple poles at $s=1/p$, $s=0$ and, in general, at $s=-k$ ($k=1, 2, \ldots\,$). However, if $p$ is an odd integer, the poles in $\Re (s)<0$ are at $s=-1, -3, \ldots$ on account of the trivial zeros of $\zeta(s)$ at $s=-2, -4, \ldots\,$. And if $p$ is an even integer there are no poles in $\Re (s)<0$. Displacement of the integration path to the left over the poles (when $p\neq 2, 4, \ldots\,$) then formally produces
\bee\label{e23}
S_{\nu,p}(a,x)=\frac{1}{\g(1+\nu)}\bl\{\frac{1}{p}\,a^{-1/p} \g(1/p)\,{}_2F_1(\f{1}{2p}, \f{p+1}{2p};1+\nu;-x^2)-\frac{1}{2}+\Upsilon(a)\br\}
\ee
with
\[\Upsilon(a)=\sum_{k=1}^\infty \frac{(-a)^k}{k!}\,\zeta(-kp)\,{}_2F_1(-\fs k,-\fs k+\fs;1+\nu;-x^2),\]
where we have used the fact that $\zeta(0)=-\fs$. From the functional relation \cite[p.~603]{DLMF}
\bee\label{e2zeta}
\zeta(s)=2^s \pi^{s-1}\,\zeta(1-s) \g(1-s) \sin \fs\pi s,
\ee
we find
\[\zeta(-kp)=-\frac{(2\pi)^{-kp}}{\pi}\,\zeta(1+kp) \g(1+kp) \sin (\fs\pi kp)\]
so that the residue sum $\Upsilon(a)$ can be written alternatively as
\[\Upsilon(a)=-\frac{1}{\pi}\sum_{k=1}^\infty \frac{(-)^k}{k!} \bl(\frac{a}{(2\pi)^p}\br)^{\!k} \zeta(1+kp) \g(1+kp) \sin (\fs\pi kp)\hspace{3cm}\]
\bee\label{e24}
\hspace{5cm}\times\ {}_2F_1(-\fs k,-\fs k+\fs;1+\nu;-x^2).
\ee

Throughout the paper we define the acute angle $\phi:=\arctan x$.
To discuss the convergence of the series (\ref{e24}) we require the large-$k$ behaviour of the above hypergeometric function.  From (\ref{a2}), this is given by
\[\frac{1}{\g(1+\nu)}\,{}_2F_1(-\fs k, -\fs k+\fs;1+\nu;-x^2)\hspace{8cm}\]
\[\sim \frac{(1+x^2)^{\frac{1}{2}k+\nu/2+3/4}}{\sqrt{\pi} (xk)^{\nu+1/2}}\,\sin ((k+\nu+\f{3}{2})\phi-\fs\pi\nu+\f{1}{4}\pi)\qquad (k\to\infty).\]
Hence we see that the late terms in $\Upsilon(a)$ are controlled in absolute value by
\[\frac{(1+x^2)^{k/2}}{k^{\nu+1/2}}\,\frac{\g(1+kp)}{k!}\,\bl(\frac{a\sqrt{1+x^2}}{(2\pi)^p}\br)^{\!k}\]
since $\zeta(1+kp) \simeq 1$ for $k\to\infty$. Consequently the residue sum $\Upsilon(a)$ converges absolutely when $p<1$ and diverges when $p>1$. When $p=1$, the sum (where only odd values of $k$ contribute) is absolutely convergent provided $a<2\pi/\sqrt{1+x^2}$ and divergent otherwise. This result may be summarised in the following theorem:
\newtheorem{theorem}{Theorem}
\begin{theorem}$\!\!\!.$ \ For $p>0$ $($when $p\neq 2, 4, \ldots$$)$, $a>0$, $x>0$, and $\nu$ real the following expansion holds:
\[S_{\nu,p}(a,x)=\frac{1}{\g(1+\nu)}\bl\{\frac{1}{p}\,a^{-1/p} \g(1/p)\,{}_2F_1(\f{1}{2p}, \f{p+1}{2p};1+\nu;-x^2)-\frac{1}{2}\]
\bee\label{e2th}
-\frac{1}{\pi}\sum_{k=1}^\infty \frac{(-)^k}{k!} \bl(\frac{a}{(2\pi)^p}\br)^{\!k} \zeta(1+kp) \g(1+kp) \sin (\fs\pi kp)\ {}_2F_1(-\fs k,-\fs k+\fs;1+\nu;-x^2).
\ee
When $p<1$, the expansion (\ref{e2th}) is an equality, but is asymptotic when $p>1$. When $p=1$, absolute convergence of the infinite series holds when $a<2\pi/\sqrt{1+x^2}$, otherwise it is divergent.
\end{theorem}
In the following subsection we analyse the case $p=1$ more carefully.
\vspace{0.3cm}

\noindent{\bf 2.1\ \ The case $p=1$}
\vspace{0.2cm}

\noindent When $p=1$, the poles in $\Re (s)<0$ are situated at $s=-1, -3, \ldots\ $. We displace the integration path in (\ref{e22}) to the left to coincide with the path $s=-2N+it$, $t\in(-\infty, \infty)$, where $N$ is a positive integer. We find
\[S_{\nu,1}(a,x)=\frac{1}{\g(1+\nu)} \bl\{\frac{1}{a}\,{}_2F_1(\fs,1;1+\nu;-x^2)-\frac{1}{2}\hspace{6cm}\]
\[\hspace{2cm}+\frac{1}{\pi}\sum_{k=0}^{N-1}(-)^k \zeta(2k+2) \,{}_2F_1(-k,-k-\fs;1+\nu;-x^2) \bl(\frac{a}{2\pi}\br)^{\!2k+1} + R_N(a)\br\},\]
where the remainder $R_N(a)$ is
\[R_N(a)=\frac{1}{2\pi i} \int_{c-2N-\infty i}^{c-2N+\infty i} \g(s) \zeta(s) \,{}_2F_1(\fs s,\fs s+\fs;1+\nu;-x^2)\,a^{-s}ds\]
\[=\frac{1}{2\pi} \int_{-\infty}^\infty \g(-2N+it) \zeta(-2N+it)\,{}_2F_1(-N\!+\!\fs it,-N\!+\!\fs\!+\!\fs it;1\!+\!\nu;-x^2) \,a^{2N-it} dt.\]
Upon observing from (\ref{e2zeta}) that
\[\g(-s)\zeta(-s)=-\frac{(2\pi)^s}{\pi}\,\frac{\zeta(1+s)}{2\cos \fs\pi s},\]
we find that
\bee\label{e2R}
|R_N(a)|<\frac{1}{4\pi} \bl(\frac{a}{2\pi}\br)^{\!2N} \zeta(1+2N) \int_{-\infty}^\infty \frac{|{}_2F_1(-N\!+\!\fs it,-N\!+\!\fs\!+\!\fs it;1\!+\!\nu;-x^2)|}{\cosh \fs\pi t}\,dt.
\ee
From (\ref{a5}) the modulus of the hypergeometric function satisfies the bound $K (1+x^2)^N e^{-\phi t}/(N^2+\f{1}{4}t^2)^{\nu/2+1/4}$, where $K$ is a positive constant and $0<\phi<\fs\pi$. Hence
\[|R_N(a)|<\frac{K}{2\pi}\bl(\frac{a\sqrt{1+x^2}}{2\pi}\br)^{\!2N} \zeta(1+2N) \int_0^\infty \frac{\cosh \phi t}{\cosh \fs\pi t}\,\frac{dt}{(N^2\!+\!\f{1}{4}t^2)^{\nu/2+1/4}}\]
\[=O\bl(\bl(\frac{a\sqrt{1+x^2}}{2\pi}\br)^{\!2N}\br).\]
Hence, as $N\to\infty$ the remainder $R_N(a)\to0$ when $a<2\pi/\sqrt{1+x^2}$. If this last condition is not met the series is asymptotic in character.
Then we have the exact result
\[S_{\nu,1}(a,x)=\frac{1}{\g(1+\nu)} \bl\{\frac{1}{a}\,{}_2F_1(\fs,1;1+\nu;-x^2)-\frac{1}{2}\hspace{8cm}\]
\bee\label{e25}
+\frac{1}{\pi}\sum_{k=0}^\infty (-)^k \zeta(2k+2) \,{}_2F_1(-k,-k\!-\!\fs;1\!+\!\nu;-x^2) \bl(\frac{a}{2\pi}\br)^{\!2k+1}\br\}
\qquad(a<2\pi/\sqrt{1+x^2}).
\ee

As an example, the special case $\nu=-\fs$, where $J_{-1/2}(anx)/(\fs anx)^{-1/2}=\cos(anx)/\sqrt{\pi}$ and
\[{}_2F_1(-k,-k\!-\!\fs;\fs;-x^2)=(1+x^2)^{k+1/2} \cos(2(k+\fs)\phi),\]
yields
\[S_{-\frac{1}{2},1}(a,x)=\frac{1}{\sqrt{\pi}}\bl\{\frac{1}{a(1+x^2)}-\frac{1}{2}+\frac{1}{\pi}\sum_{k=0}^\infty (-)^k \zeta(2k+2)\cos (2(k+\fs)\phi)\,X^{2k+1}\br\},\]
where $X:=a\sqrt{1+x^2}/(2\pi)$. Evaluation of the sum then produces
\[S_{-\frac{1}{2},1}(a,x)=\frac{1}{\sqrt{\pi}}\bl\{\frac{1}{4}\{\coth (\pi Xe^{i\phi})+\coth (\pi Xe^{-i\phi})\}-\frac{1}{2}\br\}=\frac{1}{2\sqrt{\pi}}\bl\{\frac{\sinh a}{\cosh a-\cos ax}-1\br\}\]
\[=\frac{1}{\sqrt{\pi}}\bl\{\frac{e^a \cos ax-1}{1-2e^a \cos ax+e^{2a}}\br\}.\]
This last result is readily verified to be the case by some straightforward algebra applied to the sum $\pi^{-1/2}\sum_{n\geq 1} e^{-an} \cos (anx)$.
\vspace{0.6cm}

\begin{center}
{\bf 3.  The case $p=2$}
\end{center}
\setcounter{section}{3}
\setcounter{equation}{0}
\renewcommand{\theequation}{\arabic{section}.\arabic{equation}}
In the case $p=2$, we have from (\ref{e22}) the result
\[S_{\nu,2}(a,x)=\frac{1}{\g(1+\nu)}\bl\{\frac{1}{2}\sqrt{\frac{\pi}{a}}\,{}_2F_1(\f{1}{4},\f{3}{4};1+\nu;-x^2)-\frac{1}{2}\hspace{3cm}\]
\bee\label{e30}
\hspace{3cm}+\frac{1}{2\pi i}\int_{c-\infty i}^{c+\infty i} \g(-s) \zeta(-2s)\,{}_2F_1(-\fs s,-\fs s\!+\!\fs;1\!+\!\nu;-x^2) a^sds\br\},
\ee
where in the integral we have put $s\to -s$ and $c>0$. Since there are no poles of the integrand in $\Re (s)>0$ (due to the trivial zeros of $\zeta(-2s)$) the expansion $\Upsilon(a)$ vanishes. This indicates that this contribution is exponentially small in the limit $a\to0$.

Noting from (\ref{e2zeta}) that
\[\g(-s) \zeta(-2s)=\pi^{-2s-1/2} \zeta(1+2s) \g(s+\fs),\]
we express the integral (with the further change of variable $s\to u-\fs$) as
\[\sqrt{\frac{\pi}{a}}\,\frac{1}{2\pi i}\int_{c-\infty i}^{c+\infty i} \bl(\frac{\pi^2}{a}\br)^{\!\!-u} \g(u) \zeta(2u)\,{}_2F_1(-\fs u\!+\!\f{1}{4},-\fs u\!+\!\f{3}{4};1\!+\!\nu;-x^2)\,du \qquad(c>\fs)\]
\bee\label{e31}
=\sqrt{\frac{\pi}{a}}\,\sum_{n=1}^\infty \frac{1}{2\pi i}\int_{c-\infty i}^{c+\infty i} \bl(\frac{\pi^2n^2}{a}\br)^{\!\!-u} \g(u)\,{}_2F_1(-\fs u\!+\!\f{1}{4},-\fs u\!+\!\f{3}{4};1\!+\!\nu;-x^2)\,du.
\ee
From the form\footnote{In the absence of the ${}_2F_1$ function, the integral in (\ref{e31}) can be evaluated by the Cahen-Mellin integral as $\exp\,[-\pi^2n^2/a]$; see, for example, \cite[p.~89]{PK}.} of the integrand in (\ref{e31}) we may expect a generalised Poisson-Jacobi transformation to hold for $S_{\nu,2}(a,x)$.

To demonstrate that this is the case in a particular example we consider $\nu=-\fs$, where
\[{}_2F_1(-\fs u\!+\!\f{1}{4},-\fs u\!+\!\f{3}{4};\fs;-x^2)=(1+x^2)^{u/2-1/4} \cos (u-\fs)\phi.\]
Then the integral in (\ref{e31}) becomes
\[\frac{1}{(1+x^2)^{1/4}} \frac{1}{2\pi i}\int_{c-\infty i}^{c+\infty i}
Y^{-u} \g(u)\,\cos ((u-\fs)\phi)\,du,\qquad Y:=\frac{\pi^2n^2}{a\sqrt{1+x^2}}\]
\[=\frac{1}{(1+x^2)^{1/4}} \sum_{k=0}^\infty \frac{(-Y)^k}{k!}\,\cos ((k+\fs)\phi)
=\frac{1}{(1+x^2)^{1/4}}\, e^{-Y \cos \phi} \cos (Y\!\sin \phi-\fs\phi).\]
Since
\[{}_2F_1(\f{1}{4},\f{3}{4};\fs;-x^2)=\frac{\cos \fs\phi}{(1+x^2)^{1/4}},\]
we finally have the expansion
\[S_{-1/2,2}(a,x)=\frac{1}{\sqrt{\pi}} \sum_{n=1}^\infty e^{-an^2} \cos (an^2x)=\frac{1}{\sqrt{\pi}}\bl\{\frac{1}{2} \sqrt{\frac{\pi}{a}}\,\frac{\cos \fs\phi}{(1+x^2)^{1/4}}-\frac{1}{2}\hspace{7cm}\]
\bee\label{e32}
\hspace{2cm}+\sqrt{\frac{\pi}{a}}\,\frac{1}{(1+x^2)^{1/4}} \sum_{n=1}^\infty \exp\,\bl[\frac{-\pi^2n^2}{a(1+x^2)}\br] \cos \bl(\frac{\pi^2n^2 x}{a(1+x^2)}-\frac{1}{2}\phi\br)\br\},
\ee
which is a generalised Poisson-Jacobi transformation.
This result, however, is easily verified upon use of (\ref{e12}) applied to the sum
$\pi^{-1/2} \sum_{n\geq1}e^{-an^2} \cos (an^2x)$.

\vspace{0.3cm}

\noindent{\bf 3.1\ \ The general case}
\vspace{0.2cm}

\noindent We study the general case of the integral in (\ref{e31}), viz.
\bee\label{e33}
\frac{1}{2\pi i}\int_{c-\infty i}^{c+\infty i} \bl(\frac{\pi^2n^2}{a}\br)^{\!\!-u} \g(u)\,{}_2F_1(-\fs u\!+\!\f{1}{4},-\fs u\!+\!\f{3}{4};1\!+\!\nu;-x^2)\,du\qquad (c>0)
\ee
\[=\sum_{k=0}^\infty \frac{(-\chi)^k}{k!}\,{}_2F_1(\fs k\!+\!\f{1}{4},\fs k\!+\!\f{3}{4};1\!+\!\nu;-x^2),\qquad \chi:=\frac{\pi^2 n^2}{a}\]
upon displacement of the integration path to the left over the poles of $\g(u)$. When $x^2<1$, we can series expand
the hypergeometric function to obtain
\[\sum_{k=0}^\infty\frac{(-\chi)^k}{k!} \sum_{r=0}^\infty \frac{(k+\fs)_{2r}}{(1+\nu)_r r!}\,(-x^2/4)^r
=\frac{1}{\sqrt{\pi}}\sum_{r=0}^\infty \frac{\g(2r+\fs)}{(1+\nu)_r r!} (-x^2/4)^r \sum_{k=0}^\infty \frac{(-\chi)^k (2r+\fs)_k}{k! (\fs)_k}.\]

The inner sum can be expressed as a confluent hypergeometric function in the form
\[{}_1F_1(2r+\fs;\fs;-\chi)=e^{-\chi} \,{}_1F_1(-2r;\fs;\chi)=e^{-\chi}\,\frac{(2r)!}{(4r)!}\,H_{4r}(\sqrt{\chi})\]
by application of Kummer's theorem and $H_n(x)$ is the Hermite polynomial \cite[p.~328]{DLMF}.
Thus the integral in (\ref{e33}) can be evaluated as
\bee\label{e34}
\sqrt{\frac{\pi}{a}} \sum_{n=1}^\infty e^{-\pi^2n^2/a} \,P_\nu\bl(x,\frac{\pi^2n^2}{a}\br),
\ee
where 
\[P_\nu\bl(x,\frac{\pi^2n^2}{a}\br):=\sum_{r=0}^\infty \frac{H_{4r}(\pi n/\sqrt{a})}{(1+\nu)_r r!} \,(-x^2\!/64)^r.\]
The form (\ref{e34}) (with $x^2<1$) demonstrates the Poisson-Jacobi-type structure for $S_{\nu,2}(a,x)$ in (\ref{e30}). However, we have been unable to express $P_\nu(x,\pi^2n^2/a)$ in a simpler, more recognisable form. 
\vspace{0.6cm}

\begin{center}
{\bf 4.  Two generalisations}
\end{center}
\setcounter{section}{4}
\setcounter{equation}{0}
\renewcommand{\theequation}{\arabic{section}.\arabic{equation}}
An extension of the sum in (\ref{e11}) is given by
\bee\label{e41}
S_{\nu,p}^\mu(a,x)=\sum_{n=1}^\infty e^{-an^p} \frac{J_\nu(an^px)}{(\fs an^p x)^{\nu-\mu}},
\ee
where $\mu$ is real. The same procedure employed in Section 2 yields the integral representation
\[S_{\nu,p}^\mu(a,x)=\frac{(\fs x)^\mu}{\g(1+\nu)}\,\frac{1}{2\pi i}\int_{c-\infty i}^{c+\infty i} 
\g(s+\mu)\zeta(sp)
\,{}_2F_1(\f{s+\mu}{2},\f{s+\mu+1}{2};1+\nu;-x^2) a^{-s} ds,\]
where $c>\max \{1/p,-\mu\}$. Poles of the integrand are situated at $s=1/p$ and $s=-\mu-k$, although some poles can be deleted on account of the trivial zeros of $\zeta(sp)$ and it is possible to have a double pole at $s=1/p$ when $\mu=-1/p, -1/p-1, \ldots\ $. As an example, we display the
particular case $p=2$ and $\mu\neq -\fs, -\f{3}{2}, \ldots\,$, to obtain
\[S_{\nu,p}^\mu(a,x)=\frac{(\fs x)^\mu}{\g(1+\nu)}\bl\{\frac{\g(\mu+\fs)}{2a^{1/2}}\,{}_2F_1(\fs\mu+\f{1}{4},\fs\mu+\f{3}{4};1+\nu;-x^2)\]
\bee\label{e42}
-\frac{\sin \pi\mu}{\pi} \sum_{k=0}^\infty\frac{1}{k!} \bl(\frac{a}{(2\pi)^2}\br)^{\!\mu+k}\zeta(1\!+\!2\mu\!+\!2k)\g(1\!+\!2\mu\!+\!2k)\,{}_2F_1(-\fs k,-\fs k\!+\!\fs;1\!+\!\nu;-x^2)\br\}, 
\ee
where we have made use of the functional relation in (\ref{e2zeta}). 

When $\mu=1$ it is seen that the infinite sum of residues in (\ref{e42}) vanishes, which provides us with another example of a generalised Poisson-Jacobi transformation. When $\nu=\fs$, we obtain
\[S_{1/2,2}^1(a,x)=\frac{1}{\sqrt{\pi}} \sum_{n=1}^\infty e^{-an^2} \sin (an^2x)=\frac{x}{\sqrt{\pi}}\bl\{\frac{1}{4}\sqrt{\frac{\pi}{a}}\,{}_2F_1(\f{3}{4},\f{5}{4};\f{3}{2};-x^2)+R(a)\br\}
\]
\bee\label{e43}
=\frac{1}{\sqrt{\pi}}\bl\{\frac{1}{2}\sqrt{\frac{\pi}{a}}\,\frac{\sin \fs\phi}{(1+x^2)^{1/4}}+xR(a)\br\}.
\ee
Here 
\[R(a)=\frac{1}{2\pi i}\int_{c-\infty i}^{c+\infty i} \g(1-s) \zeta(-2s) \,{}_2F_1(\f{1-s}{2},\f{2-s}{2};\f{3}{2};-x^2) a^s ds\qquad (c>0)\]
\[=\sqrt{\frac{\pi}{a}} \,\frac{1}{2\pi i}\int_{c-\infty i}^{c+\infty i} \bl(\frac{\pi^2}{a}\br)^{\!-u} \g(u) \zeta(2u) (u-\fs)\,{}_2F_1(-\fs u\!+\!\f{3}{4},-\fs u\!+\!\f{5}{4};\f{3}{2};-x^2)\,du,\]
since from (\ref{e2zeta})
\[\g(1-s)\zeta(-2s)=-\frac{(u-\fs)}{\pi^{2u-1/2}}\,\g(u)\zeta(2u)\qquad (s\to u-\fs).\]
Making use of the identity
\[{}_2F_1(-\fs u\!+\!\f{3}{4},-\fs u\!+\!\f{5}{4};\f{3}{2};-x^2)=(1+x^2)^{u/2-1/4}\,\frac{\sin ((u-\fs)\phi)}{x(u-\fs)},\]
we then obtain
\[xR(a)=-\sqrt{\frac{\pi}{a}}\,\frac{1}{(1+x^2)^{1/4}}\,\sum_{n=1}^\infty\,\frac{1}{2\pi i}\int_{c-\infty i}^{c+\infty i} \bl(\frac{\pi^2n^2}{a\sqrt{1+x^2}}\br)^{\!-u} \g(u) \sin((u-\fs)\phi)\,du\]
\[\hspace{0.2cm}=\,\sqrt{\frac{\pi}{a}}\,\frac{1}{(1+x^2)^{1/4}}\,\sum_{n=1}^\infty \sum_{k=0}^\infty \frac{(-Y)^k}{k!} \sin ((k+\fs)\phi),
\qquad Y:=\frac{\pi^2n^2}{a\sqrt{1+x^2}}.\]
Then since
\[\sum_{k=0}^\infty \frac{(-Y)^k}{k!} \sin ((k+\fs)\phi)=-e^{-Y \cos \phi} \sin(Y\!\sin \phi-\fs\phi),\]
we finally obtain
\bee\label{e44}
xR(a)=-\sqrt{\frac{\pi}{a}}\,\frac{1}{(1+x^2)^{1/4}} \sum_{n=1}^\infty \exp\,\bl[\frac{-\pi^2n^2}{a(1+x^2)}\br] \sin \bl(\frac{\pi^2n^2 x}{a(1+x^2)}-\frac{1}{2}\phi\br).\hspace{0.6cm}
\ee
Combination of (\ref{e43}) and (\ref{e44}) then gives the expansion of $S_{1/2,2}^1(a,x)$ for $a>0$ and clearly has the form of a generalised Poisson-Jacobi transformation.
This result is easily verified upon use of (\ref{e12}) applied to the sum
$\pi^{-1/2} \sum_{n\geq1}e^{-an^2} \sin (an^2x)$.

Another extension is the sum
\bee\label{e45}
T_{\nu,p}(a,x)=\sum_{n=1}^\infty e^{-an^p} \frac{I_\nu(an^p x)}{(\fs an^p x)^\nu}, \qquad x\in (0,1),
\ee
\[=\frac{1}{\g(1+\nu)}\,\frac{1}{2\pi i}\int_{c-\infty i}^{c+\infty i} \g(s)\,\zeta(sp) \,{}_2F_1(\fs s,\fs s\!+\!\fs;1\!+\!\nu;x^2)\,a^{-s}ds\qquad (c>1/p),\]
where $I_\nu(z)$ is the modified Bessel function of the first kind. Then we obtain the same expansion given in (\ref{e23}) and (\ref{e24}) with the argument of the hypergeometric function replaced by $x^2$. 
When $0<p<1$ the infinite series of residues is convergent but divergent (asymptotic) when $p>1$.

The special case $p=1$ yields
\[T_{\nu,1}(a,x)=\frac{1}{\g(1+\nu)}\bl\{\frac{1}{a}\,{}_2F_1(\fs,1;1+\nu;x^2)-\frac{1}{2}\hspace{4cm}\]
\bee\label{e46}
\hspace{4cm}+\frac{1}{\pi}\sum_{k=1}^\infty (-)^k \bl(\frac{a}{2\pi}\br)^{\!2k+1} \zeta(2k+2)\,{}_2F_1(-k,-k+\fs;1+\nu;x^2)\br\}.
\ee
From (\ref{a4}), the leading large-$k$ behaviour of the above hypergeometric function is given by
\[{}_2F_1(-k,-k+\fs;1+\nu;x^2)\sim \frac{\g(1+\nu)}{2\sqrt{\pi}}\,\frac{(1+x)^{2k+\nu+3/2}}{(kx)^{\nu+1/2}}\qquad (k\to\infty)\]
so that the infinite sum in (\ref{e46}) converges when $a<2\pi/(1+x)$. If $a(1+x)/(2\pi)=1$, we have convergence of the sum 
for $\nu>-\fs$, since it is easily seen that the hypergeometric series is positive for $k\geq 1$.
\vspace{0.6cm}

\begin{center}
{\bf 5. Concluding remarks}
\end{center}
\setcounter{section}{5}
\setcounter{equation}{0}
\renewcommand{\theequation}{\arabic{section}.\arabic{equation}}
In Theorem 1 we have presented the expansion of $S_{\nu,p}(a,x)$ for $p>0$ ($p\neq 2, 4,\ldots$) with the parameters 
$a>0$ and $x>0$. This expansion is convergent when $p<1$ but is asymptotic when $p>1$. In the case $p=1$ convergence requires the condition $a<2\pi/\sqrt{1+x^2}$. When $p$ is an even integer the character of the expansion changes
to become exponentially small. In two cases when $p=2$ it was possible to display a generalised Poisson-Jacobi transformation, although both follow in a straightforward manner from (\ref{e12}). An attempt at the general case
managed to establish the standard infinite sum of exponentials of the form $\exp\,[-\pi^2n^2/a]$ in (\ref{e34}), although the factor $P_\nu(x, \pi^2n^2/a)$ could not be simplified.

The situation when $p=4$, or higher, is more complicated and preliminary investigation suggests a composite of expansions each containing the exponential factor $\exp\,[-\pi^2n^2/a]$. We do not consider this situation further here, nor how the expansion might change as $p\to2$.

\vspace{0.6cm}

\begin{center}
{\bf Appendix: The large-$k$ behaviour of some hypergeometric functions }
\end{center}
\setcounter{section}{1}
\setcounter{equation}{0}
\renewcommand{\theequation}{\Alph{section}.\arabic{equation}}
In this appendix we determine the large-$k$ asymptotic behaviour of the hypergeometric functions appearing in the main body of the paper. We first consider the function
\bee\label{a1}
F\equiv{}_2F_1(-k,-k-\fs;1+\nu;-x^2)\qquad (k\to+\infty),
\ee
where $x^2>0$.
From \cite[pp.~388, 390]{DLMF} we have
\[F=(1+x^2)^{2k+\nu+3/2}\,{}_2F_1(k+\nu+1;k+\nu+\f{3}{2};1+\nu;-x^2)\]
\bee\label{a2}
=(1+x^2)^{2k+\nu+3/2}\,\frac{\g(1+\nu)\g(k+1)}{2\pi i\,\g(k+\nu+1)}\,\int_0^{(1+)} \frac{t^{k+\nu} (t-1)^{-k-1}}{(1+x^2 t)^{k+\nu+3/2}}\,dt,
\ee
where the integration path is a closed contour starting and finishing at $t=0$ that encircles $t=1$ in the positive sense. The above integral can be written as
\[\frac{1}{2\pi i} \int_0^{(1+)} e^{k\psi(t)} f(t)\,dt,\]
where
\[\psi(t)=\log \bl(\frac{t}{(t-1)(1+x^2t)}\br),\qquad f(t)=\frac{t^\nu}{(t-1) (1+x^2t)^{\nu+3/2}}.\]
Saddle points of $\psi(t)$ arise where $\psi'(t)=0$; that is, at $t_s=\pm i/x$. For the saddle at $t_s=i/x$, we have
\[\psi''(t_s)=\frac{2ix^3}{(1+ix)^2},\qquad f(t_s)=(i/x)^{\nu-1} (1+ix)^{-\nu-5/2},\]
\bee\label{a10}
 e^{k\psi(t_s)}=\frac{1}{(1+ix)^{2k}}=\frac{e^{-2ik\phi}}{(1+x^2)^k},\qquad \phi:=\arctan x.
\ee
The path of steepest descent through the saddle $t_s=i/x$ emanates from the origin and passes to infinity in the upper half-plane in the direction $\arg\,t=2\phi$; that is, in
$\Re (t)>0$ when $0<x<1$, along $\Re (t)=0$ when $x=1$ and in $\Re (t)<0$ when $x>1$.
The direction of integration at the saddle is $\fs\pi-\fs \arg \psi''(t_s)=\f{1}{4}\pi+\phi$; the steepest descent path through the saddle at $-i/x$ in the lower half-plane is the reflection of the path in the upper half-plane. The integration path can then be expanded to infinity around an infinite arc to pass over both saddles, which contribute equally to the integral.

Application of the saddle-point method \cite[p.~47]{DLMF} shows that the contribution to the integral from the saddle $t_s=i/x$ is
\[-\sqrt{\frac{2\pi}{k(-\psi''(t_s))}}\,f(t_s) e^{k\psi(t_s)}=- \sqrt{\frac{\pi}{k}}\,\frac{(1+x^2)^{-k-\nu/2-3/4}}{x^{\nu+1/2}}\,\exp \bl[-i(2k+\nu+\f{3}{2})\phi+\fs\pi i(\nu-\fs)\br]\]
with the conjugate expression from the saddle $t_s=-i/x$. Hence, noting that the ratio of gamma functions multiplying
the integral in (\ref{a2}) is $\g(1+\nu) k^{-\nu}$ for large $k$, we obtain the final result
\[{}_2F_1(-k,-k-\fs;1+\nu;-x^2)\hspace{8cm}\]
\bee\label{a3}
\hspace{2cm}\sim\frac{\g(1+\nu)}{\sqrt{\pi}}\,\frac{(1+x^2)^{k+\nu/2+3/4}}{(xk)^{\nu+1/2}}\,\sin \bl((2k+\nu+\f{3}{2})\phi-\fs\pi\nu+\f{1}{4}\pi\br)\qquad (k\to\infty).
\ee
In the special case $\nu=-\fs$, (\ref{a3}) yields
\[{}_2F_1(-k,-k-\fs;\fs;-x^2)\sim (1+x^2)^{k+1/2} \cos ((2k+1)\phi),\]
which is the exact result.

For the function with positive argument
\[{}_2F_1(-k,-k-\fs;1+\nu;x^2) \qquad x\in (0,1),\]
the procedure is similar, where now $\psi(t)=\log [(t/((t-1)(1-x^2t)]$ with saddles at $\pm1/x$.
The integration path can be deformed to pass over the saddle at $t_s=1/x$, where the direction of integration at the saddle is $\fs\pi$ and
\[\psi''(t_s)=\frac{2x^3}{(1-x)^2},\qquad f(t_s)=x^{1-\nu} (1-x)^{\nu-1/2},\qquad e^{k\psi(t_s)}=(1-x)^{-2k}.\]
Then we obtain the leading large-$k$ behaviour given by
\bee\label{a4}
{}_2F_1(-k,-k-\fs;1+\nu;x^2)\sim \frac{\g(1+\nu)}{2\sqrt{\pi}}\,\frac{(1+x)^{2k+\nu+3/2}}{(xk)^{\nu+1/2}}\qquad (k\to\infty).
\ee

Finally, we require an estimate for
\[|{}_2F_1(-N+\fs it,-N+\fs+\fs it;1+\nu;-x^2)|,\qquad t\in(-\infty,\infty)\] 
that appears in the estimation of the remainder term $R_N(a)$ in (\ref{e2R}). The procedure follows that employed in the estimation of (\ref{a1}) with $k$ replaced by $N-\fs it$. The modulus of the exponential factor in (\ref{a10}) becomes $e^{-\phi t}/(1+x^2)^N$, with the result that 
\bee\label{a5}
|{}_2F_1(-N+\fs it,-N+\fs+\fs it;1+\nu;-x^2)|<\frac{K (1+x^2)^N e^{-\phi t}}{(N^2+\f{1}{4}t^2)^{\nu/2+1/4}},
\ee
where $K$ is a positive constant.

\end{document}